\documentclass[12pt]{article}
\usepackage{modif}
\usepackage{amsthm,
 }
\usepackage{graphicx}
\usepackage[colorlinks=true,citecolor=black,linkcolor=black,urlcolor=blue]{hyperref}

\theoremstyle{plain}
\newtheorem{theorem}{Theorem}
\newtheorem{lemma}{Lemma}[section]
\newtheorem{corollary}[theorem]{Corollary}

\theoremstyle{definition}

\theoremstyle{remark}
\newtheorem{remark}[theorem]{Remark}

\def\N{{\bf N}}
\def\R{{\bf R}}
\def\implies{\; \Longrightarrow \;}
\let\0=\emptyset
\def\T#1,#2 {\ensuremath{T(#1;#2)}}
\def\c#1{\ensuremath{C_{#1}}}
\def\cle#1{\c{\le #1}}
\def\p#1{\ensuremath{P_{#1}}}
\def\3#1#2#3{_{#1=#2}^{#3}}   
\def\ca#1{\left|{#1}\right|}
\def\mval#1 {\delta\mathopen({#1}\mathclose)}
\def\mxval#1 {\Delta\mathopen({#1}\mathclose)}
\def\dist{\mathop{\rm dist}\nolimits}
\def\B#1#2{\ensuremath{B({#2};{#1})}}
\def\girth#1{\mathop{\rm girth}({#1})}

 \newbox\esbox 
\setbox\esbox=\hbox{\bf\raise2.1pt\hbox{.}\rm
 \hskip-1.6pt \raise.4pt\hbox{-}\hskip-2.5pt
 \bf\raise2.1pt\hbox{.}}
\def\edgesymbol{\copy\esbox}
\def\2#1,#2.{\ensuremath{{#1}\edgesymbol{#2}}}

\let\mv=\delta
\let\mx=\Delta
\def\e#1 {\ensuremath{\mathop{\rm
 e}\mathopen({#1}\mathclose)}}

\let\ep=\varepsilon
\def\br{\hfill\break}
\def\kk#1#2{\ensuremath{K_{{#1},{#2}}}}

\def\v#1{\ensuremath{\deg(#1)}}
\def\vind#1#2{\ensuremath{\deg_{#1}(#2)}}
\def\vv#1#2{\ensuremath{\deg^{\,#1}(#2)}}
\def\vvind#1#2#3{\ensuremath{\deg_{#1}^{#2}(#3)}}

\def\kdots{\ifmmode,\ldots,\else,~$\ldots\,$, \fi}
\def\row#1#2{\ifmmode #1_1\kdots #1_{#2}\else
 $#1_1$\kdots$#1_{#2}$\fi}

\def\noproof{\hfill$\square$}

\long\def\clause#1#2{\par\smallbreak\hangafter=1\hangindent
 20pt\noindent{\hbox to20pt{{#1\hfill}}{#2}}\hfill\par
 \ifdim \lastskip <\smallskipamount \removelastskip
 \penalty 15\smallskip \fi}

\def\bt{{\bf t}}
\def\bu{{\bf u}}
\def\bw{{\bf w}}
\def\bd#1{\ensuremath{{\bf d}_{#1}}}
\def\bs{{\bf s}}
\def\es{\ensuremath{\overline{\bf s}}}       
\def\ds#1{\mathop{\bf s}({#1})}
\def\ts{\ensuremath{\tilde{\bf s}}}
\def\ms#1{\mathop{\tilde{\bf s}}({#1})}
\def\V{\ensuremath{{\cal V}}}
\def\VD{\ensuremath{{\cal D}}}

\def\nb{neighbour}
\def\Tu{Tur\'an }

\title{\bf Sizes of the extremal girth five graphs
of orders from 40 to 49.}

\author{J\"orgen Backelin\\
\small Department of Mathematics\\[-0.8ex]
\small University of Stockholm\\
\small\tt joeb@math.su.se\\}

 \date{\today}

\begin{document}
 \maketitle

 \begin{abstract}
 The \Tu type numbers for graphs without 3-cycles
and 4-cycles are determined for vertex numbers from
40 to 49 inclusive. Hence, now, 43 of the first 50
numbers of OEIS~{\bf A006856} are known. Estimates
for the remaining seven numbers are presented.

 \end{abstract}

 \section{Introduction.}
 \label{S:intr}

 In this note, all considered graphs $G$ are simple
and undirected, with orders (numbers of vertices)
$n(G)$ and sizes (numbers of edges) $e(G)$,
respectively. The main object is the extremal graphs
of girth at least five, and in particular their
sizes. Thus, we consider the \Tu type numbers
 $$\T\cle4,n = \max (e \;:\; \exists\, \hbox{ graph
$G$ of order $n$, size $e$, and girth } \ge 5).$$
 They form item {\bf A006856} in \emph{The on-line
encyclopedia of integer sequences} (\cite{OEIS}),
and for $n\le32$ were determined by Garnick, Kwong,
Lazebnik, Nieuwejaar, McKay, Codish, Miller,
Prosser, and Stuckey (\cite{GN}, \cite{GKL},
\cite{CMPS}).

 Moreover, Dutton and Brigham (\cite{DB}), and
independently the authors of \cite{GKL}, also found
a ``theoretical upper bound''
 \begin{equation}
 \label{E:tub}
 e(G) \le 0.5n(G)\sqrt{n(G)-1},
 \end{equation}
 and the latter authors also noted that this bound
is attained for the rather symmetric
Hoffman-Singleton graph $HS$ (\cite{HS}); whence
they deduced that
 $$\T\cle4,50 = e(HS) = 175.$$
 Now, in my experience, the series of extremal
graphs of increasing orders for similar kinds of
conditions often contain unique and highly symmetric
graphs for some orders; and, when they do, all the
extremal graphs of the nearest preceeding orders
usually are induced subgraphs of these symmetric
graphs. Thus, I suspected the same to be true in
this case.

 In this note, I prove this, but probably not to the
fullest possible extent. The two first theorems do
imply that for orders $40 \le n(G) \le 49$ there are
extremal graphs which are subgraphs of $HS$; but I
have succeeded to prove that there also are no other
extremal graphs only for $n \in \{40,45,47,48,49\}$.
For these $n$, the extremal graphs also are unique
(up to isomorphisms). Unicity does not hold for all
the remaining 5 values, but I find it likely that
also in these cases all the extremal graphs are
subgraphs of Hoffman-Singleton graphs.
 Moreover, for several lower values of $n$, the
corresponding questions are open.

 Recall that the vertex set $V(HS)$ of $HS$ may be
partitioned into two parts $V'$ and $V''$, such that
the induced $HS$ subgraph $HS[V^{(\nu)}]$ on either
part is isomorphic to the disjoint union $5\c5$ of
five 5-cycles, and that the induced subgraphs on the
unions of one 5-cycle from each part are Petersen
graphs.
 If we choose one 5-cycle from $V''$, and
successively remove its vertices from $HS$, the
resulting induced $HS$ subgraphs realise the lower
bounds
 \begin{equation}
 \label{E:hi}
 \T\cle4,45 \ge 145, \hbox{ and }
 \T\cle4,n \ge 6n-126 \hbox{ for } n = 46\kdots49.
 \end{equation}
 (All but the first one of these lower bounds also
were found by Garnick, Kwong, and Lazebnik in
\cite{GKL}.)

 Note, that when the entire 5-cycle is removed, the
resulting graph $G$ has minimal degree $\mval G =
6$, maximal degree $\mxval G = 7$, induced graph
$G_6 = G[V_6] \simeq 5\c5$, and induced graph $G_7 =
G[V_7] \simeq 4\c5$, where $V_i = V_i(G) := \{v \in
V(G) : \deg v = i\}$. If we repeat the procedure, by
choosing one of the 5-cycles in $G_6\,$, and
successively removing its vertices, we get further
$HS$ subgraphs, yielding
 \begin{equation}
 \label{E:lo}
 \T\cle4,40 \ge 120, \hbox{ and }
 \T\cle4,n \ge 5n-81 \hbox{ for } n = 41\kdots44.
 \end{equation}

 The main object of this article is to prove
equalities in (\ref{E:hi}) and (\ref{E:lo}), i.e.,
to prove
 \begin{theorem}
 \label{T}
 $\T\cle4,45 = 145$, and
 $\T\cle4,n = 6n-126$ for $n = 46\kdots49$.
 \end{theorem}
 and
\begin{theorem}
 \label{T2}
 $\T\cle4,40 = 120$, and
 $\T\cle4,n = 5n-81$ for $n = 41\kdots44$.
 \end{theorem}

 The main technical tools for proving the theorems
are the following main lemmata:

 \begin{lemma}
 \label{L:m} 
 If $G$ is a graph with $\girth G \ge 5$, $n(G)=45$,
and $e(G)=145$, then $\mval G = 6$, $\mxval G = 7$,
and each vertex has exactly two \nb s with the same
degree as itself.
 \end{lemma}
 and
 \begin{lemma}
 \label{L:m2} 
 If $G$ is a graph with $\girth G \ge 5$, $n(G)=40$,
and $e(G)=120$, then $G$ is 6-regular (and hence a
cage).
 \end{lemma}

 A more precise structural result is
\begin{theorem}
 \label{T:HSsub}
 For any fixed order $n \in \{40,45,47,48,49\}$, the
girth $\ge5$ graphs of maximal size are isomorphic,
and subgraps of Hoffman-Singleton graphs.
 \end{theorem}
 The case $n=40$ follows directly from
lemma~\ref{L:m2}, and from P.-K. Wong's precise
characterisation of (6,5) cages as
`Hoffman-Singleton minus Petersen' graphs in \cite
W. For the higher $n$ values, the claims are
consequences of lemma~\ref{L:m} and theorem~\ref T,
as seen in section~\ref{S:unic}.

 \begin{remark}
 The unicity part of theorem~\ref{T:HSsub} may be
reformulated as
 $$a(40) = a(45) = a(47) = a(48) =
a(49) = 1,$$
 where $a$ is the function listed in the OEIS item
{\bf A159847} (\cite{OEIS2}).
 \end{remark}

 \begin{remark}
 The order in which the results are listed probably
appears to be counter-intuitive. However, I found
the analysis of \T\cle4,n \ easier for values closer
to 50, than for those a bit lower. Moreover, both
the auxiliary lemma~\ref{L:m} and theorem~\ref T are
proved by ``intrinsic'' means, which do not depend
on lower \T\cle4,n \ values or bounds, and nor on
structure results for small extremal girth~5 graphs.

 In fact, instead of investigating the rather large
number of potentially possible degree sequences one
by one, I have as far as I was able `linearised' the
influence of differences in degree sequences.
 In an intermediate step, we shall work with some
`virtual degree sequences', where actually some
entries may be negative.
 This enables a reduction of the main part of the
proof to a kind of `elementary calculus', rather
than a division into a cumbersome number of cases.
In particular, the proof presented here does not in
any manner depend on computer calculations.

 On the other hand, I have found no independent way
to prove the corresponding auxiliary result
lemma~\ref{L:m2}, but instead partly had to resort
to the traditional `recursive' use of established
values of and upper bounds for \T\cle4,n \ for
$n<40$, and to case divisions. The best bounds (to
my knowledge) for `intermediate' $n$ (i.e., for $n =
33\kdots39$) are given in table~\ref{tab:im}; they
contain a few improvements. However, it should be
noted that such bounds to some extent are based on
exact values which seem to be announced but not
proven in the literature.

 Notwithstanding, theorem~\ref{T2} does not depend
on theorem~\ref T; but the former might be used to
abbreviate the proof of a minor part
lemma~\ref{L:m}, by proving (\ref{E:fpos}) faster.
 \end{remark}

 \section{Auxiliary results.}
 \label{S:aux}

 For any graph $G$ and vertex $v \in V(G)$, let
$S(v) = \{w \in V(G) : \2v,w. \in E(G)\}$ (the `unit
sphere' centred at $v$, or the open \nb hood of
$v$), $B(v) = \B1v = S(v) \cup \{v\}$ (the `unit
ball' centred at $v$, or the closed \nb hood of
$v$), $\B2v = \bigcup\limits_{w\in S(v)} B(w)$ (the
`radius two ball' centred at $v$), and
 $$\vv2v = \vvind G2v = \sum\limits_{w\in S(v)} \v
w$$
 (the `second degree' of $v$).

 For any $W_1,W_2,W_3 \subseteq V$, let
$p_3(W_1,W_2,W_3)$ be the number of ordered paths of
order~3 (and thus size~2), which have the $i$'th
vertex in $W_i$ for $i=1,2,3$. Our main application
of this counter will be $p(U,W) := 0.5 p_3(W,U,W)$,
the number of (unordered) 3-vertex paths with middle
vertex in $U$ and end vertices in $W$ (or,
equivalently, with middle vertex in $U$ and both
edges in $E(U,W)$).

 \smallskip
 If $G$ is a graph of girth $\ge5$, and $v \in
V(G)$, then any two different paths of size two and
starting from $v$ will have different end points.
This simple observation is sufficient for deducing
the following fairly well-known lemmata.

  \begin{lemma}
 \label{L:vv}
 If $G$ is a graph with $\girth G \ge 5$ and $v \in
V(G)$, then $\ca{\B2v} = \vv2v+1$. In particular,
thus, $\vv2v \le n(G)-1$.\noproof
 \end{lemma}

 \begin{lemma}
 \label{L:delbort}
 If $G$ is a graph with $\girth G \ge 5$, $v \in
V(G)$, $l \in \N$, and \row wl are distinct $S(v)$
elements, then $e(G-\{v,\row wl\}) = e(G) - \v v -
\sum\3i1l \v{w_i} + l$. In particular, then $e(G)
\le \T\cle4,n-l-1 + \v v + \sum\3i1l
\v{w_i}-l$.\noproof
 \end{lemma}

 \begin{lemma}
 \label{L:helbort}
 If $G$ is a graph with $\girth G \ge 5$ and $v \in
V(G)$, then $e(G-B(v)) = e(G)-\vv2v$. In particular,
$e(G) \le \T\cle4,n-{\mxval G }-1 + n(G)-1$.\noproof
 \end{lemma}

 \begin{lemma}
 \label{L:stigbort}
 If $G$ is a graph with $\girth G \ge 5$ and $P$ is
a path of order $r$ (and thus size $r-1$) in $G$,
then $\sum\limits_{v\in V(P)} \v v \ge e(G) -
\T\cle4,n-r +r-1$.\noproof
 \end{lemma}

 By counting 3-vertex walks either by their first
or by their middle vertex, we get
 \begin{lemma}
 \label{L:vvsq}
 $\displaystyle\sum_{v\in W} \vvind {G[W]}2v =
p_3(W,W,W)-2e(W) = \sum_{v\in W} \vind{G[W]}v^2$,
for any\break
 $W \subseteq V$.\noproof
 \end{lemma}

 If \p3 is a 3-vertex path with vertex set $\{p_1,
p_2, p_3\}$, of which $p_1,p_3 \in W$, then $\{p_1,
p_3\}$ is one of the ${\ca W\choose2}$ 2-sets of
vertices in $W$; however, not one of the $e(W)$
edges, and nor two common \nb s of a third element
in $V$. Thus, we have
 \begin{lemma}
 \label{L:p-estim1}
 If $G$ is a graph with $\girth G \ge 5$ and $U,W
\subseteq V$, then\br
 ${\ }$\hfill $\displaystyle p(U,W) \le p(V,W) \le
{\ca W\choose2}-e(G[W])$.\noproof
 \end{lemma}

 \begin{lemma}
 \label{L:p-estim2}
  If $G$ is a graph with $\girth G \ge 5$, $U,W
\subseteq V$, and $U$ is the disjoint union of $U'$
and $U''$, then $p(U,W) = p(U',W)+p(U'',W)$.\noproof
 \end{lemma}

 We mainly shall apply $p(U,W)$ analysis in
situations where moreover $U$ and $W$ are disjoint.
In these cases, $p(U,W) \le {\ca
W\choose2}-e(W)-p(W,W)$ by lemma~\ref{L:p-estim2};
whence $p(W,W)$ estimates provide upper bounds of
$p(U,W)$. For a lower bound, given $\ca U = m$ and
$\ca{E(U,W)} = \ca {\{\hbox{edges between $U$ and }
W\}} = z$ (say), note, that $p(U,W)$ is minimal, if
the $\ca {W\cap S(u)}$ are distributed as evenly as
possible, when $u$ runs through $U$. Thus, if in
addition $rm \le z \le (r+1)m$ for an integer $q$,
then
 $$p(U,W) \ge ((q+1)m-z){q\choose2} +
(z-qm){q+1\choose2},$$
 since the RHS (right hand side) equals $p(U,W)$, if
$\ca {W\cap S(u)} \in \{r,r+1\}$ for each $u$ in
$U$. Actually, a little reflection should convince
the reader that this inequality also holds, if $z <
qm$ or $z > (q+1)m$, although then one of the RHS
terms is negative. Hence, without further
restrictions, the following lemma holds:
 \begin{lemma}
 \label{L:p-estim3}
 If $G$ is a graph with $\girth G \ge 5$, $V'$ and
$V''$ are disjoint subsets of $V(G)$, and $r$ and
$s$ are natural numbers, then
 \begin{eqnarray*}
 &((r+1)\ca U-\ca{E(U,W)}){r\choose2} +
(\ca{E(U,W)}-r\ca U){r+1\choose2}\\
 \le& p(U,W)\\
 \le& {\ca W\choose2} - e(G[W]) - ((s+1)\ca
W-2e(G[W])){s\choose2} - (2e(G[W])-s\ca
W){s+1\choose2}&{\qquad}\square
 \end{eqnarray*}
 \end{lemma}

 \section{Proof of the first main lemma.}
 \label{S:ml}

 For each potentially possible vertex degree $i$,
let $V_i = \{v \in V | \deg v = i\}$, and let $n_i =
\ca{V_i}$ (the cardinality of $V_i$). The degree
sequence is $\bs = \ds G = (n_0, n_1\kdots n_{44})$,
which may be considered as an element in $\R^{45}$.
(In practice, of course, most entries in \bs\ are
zero.)

 As a vector, \bs\ satisfies the linear restrictions
 $\sum_i n_i = 45$ and $\sum_i in_i = 2e = 290$, or,
in other words,
 \begin{equation}
 \label{E:V-cond}
 \bs\cdot\bu = 45 \hbox{ and } \bs\cdot\bw = 290,
 \end{equation}
 where $\bu = (1,1\kdots1)$ and $\bw =
(0,1,2\kdots44)$.
 Let $\V = \{\bt \in \R^{45} | \bt\cdot\bu = 45$
and $\bt\cdot\bw = 290\}$, the set of \emph{virtual
degree sequences}, and let $\VD = \{\bt \in \R^{45}
| \bt\cdot\bu = \bt\cdot\bw = 0\}$, the linear
$\R^{45}$ subspace of \emph{(degree sequence)
deviations}. Clearly, $\bt,\bt' \in \V \implies
\bt-\bt' \in \VD$, and more precisely
 \begin{equation}
 \label{E:V-class}
 \V = \es+\VD,
 \end{equation}
 where $\es = (\cdots0,25,20,0\cdots)$ is the
`globally most even degree distribution', satisfying
$\ds G = \es \iff \mval G = 6 \land \mxval G = 7$.

 Put $V' := V_{\le6} := \{v' \in V : \deg v \le 6 \}
= \bigcup_{i\le6} V_i$, $V'' := V_{\ge7} := \{v' \in
V : \deg v \ge 7 \} = \bigcup_{i\ge7} V_i$, $G'' =
G[V'']$, $n'' = \ca {V''} = n(G'')$, $n' = \ca{V'} =
45-n''$, $z = \ca {E(V',V'')}$, $e'' = e(G'')$, $e'
= e(G')$, $a = \sum_{v\in V'} 6-\deg v =
\sum_{i\le5} (6-i)n_i$, and $b = \sum_{v\in V''}
\deg(v)-7 = \sum_{i\ge8} (i-7)n_i$, and for each $v
\in V$, let $\deg''(v) $ be the number of $v$ \nb s
in $V''$.

 The idea is to give a lower estimate $p_l$ and an
upper estimate $p_u$ for the quantity $p :=
p(V',V'')$ (whence necessarily $p_l \le p \le p_u$),
and to show that, on the other hand, always $p_l \ge
p_u$, with equality only for the prescribed form of
$G$. For our first estimates, we let $p_l$ and $p_u$
be the lower and the upper $p$ bounds in
lemma~\ref{L:p-estim3}, with $r=4$ and $s=2$,
respectively. This yields
 $$p_l = 4z-10n' = 6n'+4(z-4n'),$$
 and
 $$p_u = {n''\choose2} - 14.5n'' - 2.5b + 2.5z,$$
 where for the second equality we also use the fact
that
 $$e'' = \left(\sum_{v\in V''} \v v - \ca{E(V',V'')}
\right) / 2 = 0.5n''+0.5b-0.5z\,.$$
 We thus indeed have $p_l \le p \le p_u$, where the
second inequality is an equality if and only if $2
\le \deg''(v) \le 3$ for each $v \in V''$.

 For $i \ne 6,7$, we define the elementary deviation
$\bd i = (d_{i,0}, d_{i,1}\kdots d_{i,44}) \in \VD$
in the following manner:\\
 $\bd5 = (0,0,0,0,0,1,-2,1,0\cdots)$ and $\bd8 =
(0,0,0,0,0,0,1,-2,1,0\cdots)$.\\
 For $0\le i\le4$, let $d_{i,i} = 1$, $d_{i,5} =
-(6-i)$, $d_{i,6} = 5-i$, and all other $d_{i,j} =
0$. Similarly, for $9\le i\le44$, let the
non-zero entries of \bd i be $d_{i,7} = i-8$,
$d_{i,8} = -(i-7)$, and $d_{i,i} = 1$.\\
 Then, elementary vector calculation shows that
 \begin{equation}
 \label{E:dev}
 \ds G-\es = \sum\3i04 n_i\bd i + a\bd5 + b\bd8 +
\sum\3i9{44} n_i\bd i \in \VD.
 \end{equation}
 Note, that all \bd* coefficients in the right hand
side are non-negative.

 Let $\ms G =
(0,0,0,0,0,a,n'-2a+b,n''+a-2b,b,0\cdots) = \es +
a\bd5 + b\bd8 \in \V$, the \emph{trunkated virtual
degree sequence of $G$}, and note that $\ds G=\ms G$
if and only if $5 \le \mval G $ and $\mxval G \le
8$.

 For any $\bt = (t_0\kdots t_{44}) \in \R^{45}$, we
may put $n'(\bt) = \sum\3j06 t_j$, $n''(\bt) =
\sum\3j7{44} t_j$, $a(\bt) = \sum\3j05 (6-j)t_j$,
and $b(\bt) = \sum\3j8{44} (j-7)t_j$, making $n'(G)
= n'(\ds G)$, {\sl et cetera}.
 Note, that these four functions are linear, and
have $n'(\es)=25$, $n''(\es)=20$, and $a(\es) =
b(\es) = 0$. Moreover, $n'(\bd5) = -1$, $n''(\bd5) =
1$, $a(\bd5) = 1$, and $b(\bd5) = 0$; $n'(\bd8) =
1$, $n''(\bd8) = -1$, $a(\bd8) = 0$, and $b(\bd8) =
1$; and $n'(\bd i) = n''(\bd i) = a(\bd i) = b(\bd
i) = 0$ for all other $i$.

 In particular, $n''(\bs-\es)+b(\bs-\es)-a(\bs-\es)
= 0$ for each $\bs-\es \in \VD$, whence
 \begin{equation}
 \label{E:20}
 n''(\bs)+b(\bs)-a(\bs) =
 n''(\es)+b(\es)-a(\es) = 20+0+0 = 20
 \hbox{ for each } \bs \in \V.
 \end{equation}

 Finally, also put $f = f(G) := 20-n'' = b-a$.

 \bigskip
 We now collect a few preliminary results, under
these assumptions and with this notation.

 \begin{lemma}
 \label{L:z}
 $z \ge 100+2b$.
 \end{lemma}

 \begin{proof}
 For each $v \in V''$, since $\vv2v \le 44$ by
lemma~\ref{L:vv}, there are at least $7\deg v-44$
edges \2v,w. between $v$ and $V'$, ``counted with
multiplicity'', where the multiplicity of such a
\2v,w. is $7-\deg w$. Thus and by~(\ref{E:20}),
indeed
 $$z \ge 5n''+7b-\sum\3i06 (6-i)in_i \ge 5n''+7b-5a
= 5(n''+b-a)+2b = 100+2b.$$
 \end{proof}

 Since on the other hand $z \le 6n'-a \le 6n'$, we
directly get
 \begin{corollary}
 \label{C:z}
 $n' \ge \left\lceil\frac{100}6\right\rceil = 17
\implies n'' \le 28.$\noproof
 \end{corollary}

 Note, that
 $$p_u-p_l = {n''\choose2} - 14.5n'' - 2.5b + 2.5z -
(4z -10n') = 450 + {n''\choose2} - 24.5n'' - 2.5b -
1.5z\,.$$
 For any fixed \ts, this is a decreasing function of
$z$, whence and by lemma~\ref{L:z} we get
 $$p_u-p_l \le 450 + {n''\choose2} - 24.5n'' - 2.5b
- 150 -3b = 0.5 \cdot \bigl( (n''-20)(n''-30) - 11b
\bigr);$$
 and thus have deduced
 \begin{lemma}
 \label{L:n}
 $11b \le (n''-20)(n''-30) = f^2+10f$.\noproof
 \end{lemma}

 By lemma~\ref{L:vv} $\sum\limits_{w\in V'\cap S(v)}
(6-\v w) \ge 6\v v-44$ for any $v\in V(G)$; whence
we have
 \begin{lemma}
 \label{L:a}
 $a \ge 6\mx-44$.\noproof
 \end{lemma}

 Thus, and by direct counting and~(\ref{E:20}),
 \begin{lemma}
 \label{L:mx}
 $\mx \ge 7 + \left\lceil\frac b{n''}\right\rceil =
7 + \left\lceil\frac{(20-n'')+a}{n''}\right\rceil
\ge \frac{140-6f+a}{20-f} \ge 6 +
\frac{6\mx-24}{20-f}$.\noproof
 \end{lemma}

 \bigskip
 We now may prove lemma~\ref{L:m}. First, by
corollary~\ref{C:z}, $f \ge -8 > -10$, but $0 \le
f(f+10)$ by lemma~\ref{L:n}. Hence
 \begin{equation}
 \label{E:fpos}
 f\ge0,
 \end{equation}
 and if $f=0$, then we indeed must have ``equalities
everywhere'', and may deduce the conclusions of the
lemma.

 Thus, when we, for a while and for a contradiction,
assume that $G$ satisfies the prerequisites but not
all the conclusions in lemma~\ref{L:m}, then $ 1 \le
f = 20-n'' \le 19$, and
 \begin{equation}
 \label{E:b1}
 b = f+a \ge 1+0 = 1\,.
 \end{equation}
 Thus $\mx \ge 8$ by lemma~\ref{L:mx}, and moreover
we may eliminate $a$ and $b$ from the inequalities
in lemmata~\ref{L:n} and~{L:a}, and get
 $$\frac{1+\sqrt{264\mx-1935}}2 \le f,$$
 which together with lemma~\ref{L:mx} (and the
bounds we already have deduced for $f$ and $\mx$)
yields
 $$11\mx^3-243\mx^2+1728\mx-4344 \le 0,$$
 and thus that $\mx=8$. However, then $f \ge 0.5
\lceil 1+\sqrt{177} \rceil = 8$ $a \ge 6\cdot8-44 =
4$, $f \ge 8+4 = 12$, and $\frac b{20-f} \le 1$,
whence we must have equalities `everywhere'. Thus,
in particular, $V' = V_8$ and $n_8 = f = n'' = 12$,
$V_{\le5} \ne \0$, and $\deg''(v) = 12 > 5 \ge \v v$
for each $v \in V_{\le5}$, the sought
contradiction.\noproof

 \section{Proof of theorem~\ref T.}
 \label{S:thp}

 The lower bounds are proven in the introduction.

 Now, assume for a contradiction that $n$ were
minimal in $\{45,46,47,48,49\}$ with the property
that \T \cle4,n \ were strictly larger than the
lower bound given in~(\ref{E:hi}). Choose a $G$ with
$\girth G \ge 5$, $n(G)=n$, and $e(G)$ exactly one
more than that bound. Since then $e(G) < 3.5n$,
there is a vertex $v \in V(G)$ with $\deg v = \mval
G \le 6$. Now, if $\mval G > 0$ then let $w$ be a
\nb\ of $v$ and put $u := v$; else, choose any edge
\2u,w. in $E(G)$. Put $G' := G-\2u,v.$, and note,
that therein $v$ has degree $\mval G' \le 5$.

 However, if $n=45$, then $\mval G' = 6$ by
lemma~\ref{L:m}. Thus, instead, $n\ge46$. By the
minimality of $n$, the conclusions of the theorem
would hold for $G'' := G-v$. Since on the other hand
$e(G'') = e(G)-\deg v \ge e(G)-6$, we must have
$n=46$, $e(G)=151$, and $G''$ must satisfy the
assumptions and thus the conclusions of
lemma~\ref{L:m}.

 In particular, $G''$ would contain 25 vertices of
degree 6, and each one of these would have two
neighbours of the same degree. Pick any such vertex
$v'$ and \nb\ $w''$, such that neither $v'$ nor
$w''$ were adjacent to $v$. Then \2v',w'. were an
edge between two vertices of degree 6 in $G$; whence
$G-v'$ were a graph with $\girth{G-v'} \ge 5$,
$n(G-v') = 45$, and $e(G-v')=145$, but with $\mval
G-v' = 5$, in contradiction to lemma~\ref{L:m}.
 \noproof

\section{Unicity.}
 \label{S:unic}

 We now study the more precise structures of
extremal graphs realising the \Tu type numbers we
just have established. First, again consider a $G$
with $\girth G \ge5$, $n(G) = 45$, and $e(G) = 145$.

 By lemma~\ref{L:m}, both $G_6$ and $G_7$ are
2-regular, and thus consist of disjoint unions of
cycles. We start by proving that each one of these
cycles has length~5; i.e., that $G_7 \simeq 4\c5$
and $G_6 \simeq 5\c5$; and continue by determining
$E(V_6,V_7)$, and the sets of vertices with mutual
distances 3.
 
 Let $\c\ell$ be a $G_7$ component, and, for a
contradiction, assume $\ell\ne5$, whence $\ell\ge6$.
Put $\{\row c\ell\} := V(\c\ell)$, $X_i := V_6 \cap
S(c_i)$, and $Y(v) := V(\c\ell) \cap S(v)$, for $1
\le i \le \ell$ and $v \in V_6\,$.
 For $\ell = 6$ (or $\ell=7$), each $\ca{Y(v)} \le
2$, with equality in at most 3 (7) cases, causing
 $$25 = n_6 \ge 5\ell - 3\ (5\ell-7) = 27\ (28,\
\hbox{ respectively),}$$
 in either case a contradiction; whence instead
$\ell\ge8$. Now consider $(X_{i+1}\kdots X_{i+8})$,
for any fixed $i$ (counting indices modulo $\ell$).
Then, for each higher index, we get at least 5, 5,
5, 4, 3, 2, 1, or 0 `new' elements (i.~e., elements
in the respective $\displaystyle X_{i+j} \setminus
\bigcup\3k1{j-1}X_k$), respectively; since the sum
of these amounts is $25 = n_6$, we must have
equalities. In particular, $X_{i+8}$ has one
distinct member in each one of $X_{i+1}\kdots
X_{i+5}\,$. Analogously, $X_i$ has one distinct
member in each of $X_{i+3}\kdots X_{i+8}$. In
particular, $X_1 \cap X_8 \ne \0 \implies
\ell\ge10$. Moreover, hence, each $Y(v) \ne\0$, the
`index gap lengths modulo $\ell$' in $Y(v)$ all
belong to $\{3\kdots7\}$, and
 $$\forall\,i \in \{1\kdots\ell\}\, \forall\, j \in
\{3\kdots7\} \exists!\, v \in V_6 : c_i,c_{i+j} \in
Y(v).$$

 Applying this for $i,j=3$, there were a $v \in V_6$
such that $c_3,c_6 \in Y(v)$. If $\{v',v''\} := V_6
\cap A(v)$, then $Y(v) \cap Y(v^{(\nu)}) = \0$, and,
in fact, $(c_i \in Y(v) \land c_k \in Y(v^{(\nu)})
\implies \ca{i-k} \ge 2)$, for $\nu=1,2$; forcing
$c_1,c_8 \in Y(v') \cap Y(v'')$, and contradiction.

 Thus, instead, indeed, $G_7 \simeq 4\c5$. Name the
$G_7$ components \c5, $\c5'$, $\c5''$, and $\c5'''$,
with $V(\c5^{(\nu)}) = \{\row{c^{(\nu)}}5\}$ for
$\nu=0\kdots3$. Moreover, let $X_i^{(\nu)} := V_6
\cap S(c_i^{(\nu)})$. For each $\nu$ and~$i$,
$\{X_i^{(\nu)}\}\3i15$ is a 5-partition of $V_6$,
such that each vertex in an $X_i^{(\nu)}$ has its
two $V_6$ \nb s in $X_{i-2}^{(\nu)}$ and
$X_{i+2}^{(\nu)}$. Hence, for any $G_6$ component
$\c\ell$, exactly every fifth vertex belongs to
$X_1^{(\nu)}$; whence $5 | \ell$. Moreover, $1 \ge
\ca {\c\ell\cap(X_1\cup X'_1)} = 0.2\ell$, whence
indeed $\ell=5$.

 Thus, indeed, $G_6 \simeq 5\c5$, and each
$X_i^{(\nu)}$ intersects each one of these five \c5
in exactly one vertex. In fact, for each choice of
one of the $G_7$ components and one of the $G_6$
components, the induced subgraph on these ten
vertices is a Petersen graph. This makes it natural
to present the $G_6$ components as pentagrams rather
than pentagons. Thus, we may let these components be
$D_l = \{d_{l,1}\kdots d_{l,5}\}$ for $l=0\kdots4$,
with the unusual prescriptions that their edges be
\2d_{l,i},d_{l,i+2}., with the second index
interpreted modulo 5.

 In this manner and without loss of generality, we
get
 $$X_i = \{d_{0,i},\kdots d_{4,i}\}$$
 for $i=1\kdots5$. We also may reindex the vertices
in each remaining complonent $\c5^{(\nu)}$
$(\nu=1,2,3)$ of $G_7$, in such a way, that
 $$\bigcap\3\nu03 X^{(\nu)}_i = \{d_{0,i}\},\
i=1\kdots5.$$
 As a consequence (and all the time employing the
girth condition), there must be a `shift function'
$\phi : \{1,2,3\} \times \{1,2,3,4\} \longrightarrow
\{1,2,3,4\}$, such that, for each $\nu \in
\{1,2,3\}$, $l \in \{1,2,3,4\}$, and $i \in
\{1,2,3,4,5\}$, we have
 $$X^{(\nu)}_i \cap D_l = \{d_{l,i+\phi(\nu,l)}\}.$$
 If necessary, by means of some further reindexing,
we now may determine $G$ uniquely. Indeed, for some
$l'$, both $d_{l',1}$ and $d_{l',2}$ must share
$V_7$ \nb s with $d_{0,1}$, i.~e., they must belong
to $\bigcup_\nu X^{(\nu)}_1$; and we may rearrange
the $l$ and the $\nu$ to have $l'=1$, $\phi(1,1)=1$,
and $\phi(2,l) \equiv 2l \pmod5$ for $l=1\kdots4$.
Now, since
$\{c''_1,d_{0,1},c'_1,d_{2,1+\phi(1,2)}\}$ is not
the vertex set of a \c4 in $G$, we must have
$\phi(1,2) \ne 4$; and similarly considering
 $$\{c_2,d_{1,2},c'_1,d_{2,1+\phi(1,2)}\}\ \hbox{
 and } \{c''_5,d_{1,2},c'_1,d_{2,1+\phi(1,2)}\}$$
 yields $\phi(1,2) \ne 1,3$; whence instead
$\phi(1,2) = 2$. In the same manner, we get
$\phi(1,l) = l$, for all $l$. Furthermore, $D_1 \cap
X'''_1$ is either $\{d_{1,4}\}$ or $\{d_{1,5}\}$,
i.e., $\phi(2,1) \in \{3,4\}$; and similar girth
analysis as before shows that then either $\phi(3,l)
\equiv 3l \pmod5$ or $\phi(3,l) \equiv 4l \pmod5$,
for all $l$. Finally, the second case by some
reindexing can be seen to be isomorphic to the first
one.

 Thus, we get $\phi(\nu,l) \equiv \nu l$; or, in
other words, up to isomorphisms, we have
 $$E(V_7,V_6) = \{\2c^{(\nu)}_i,d_{l,i+\nu l}. :
i=1\kdots5,\; \nu=0\kdots3,\; l=0\kdots4\};$$
 which completely determines $G$. Moreover, the only
pairs of elements in $G$ of distance greater than 2
are the pairs belonging to the same part in a
5-partition of $V_6$, consisting of ``the missing
$X_*^{(*)}$'', namely $\row {X^{(4)}}5$, where
 $$d_{l,k} \in X_i^{(4)} \iff k \equiv i+4l
\pmod5.$$

 Now, it is obvious that $G$ can be extended to a
Hoffman-Singleton graph, by adding a new \c5 with
vertices $\row{c^{(4)}}5$, and with further edges
from $c_i^{(4)}$ to the $X_i^{(4)}$ elements.
However, we may do better. If $G'$ is a girth 5
supergraph of $G$ with $45 < n(G') < 50$ and $e(G')
= \T\cle4,n(G') $, then each $v \in V(G') \setminus
V(G)$ must have several \nb s in $V(G)$; and these
\nb s must have distance 3 in $G$ and thus belong to
the same $X_i^{(4)}$. This yields
 \begin{lemma}
 \label{L:sup-HS}
 If $G$ is a graph with $\girth G \ge 5$, $45 \le
n(G) \le 50$, $e(G) = \T\cle4,n(G) $, and containing
a subgraph of order 45 and size 145, then $G$ is a
subgraph of a Hoffman-Singleton graph.\noproof
 \end{lemma}

 This includes the $n=45$ case of
theorem~\ref{T:HSsub}. We also get the cases
$n=47,48,49$, as soon as we prove that extremal
graphs of these orders indeed contain extremal
subgraphs of order 45. (Note, however, that the
corresponding statement for $n=46$ is not true;
removing a bipartite graph \kk13 from $HS$ yields an
extremal order 46 graph, all of whose order 45
subgraphs have sizes $<145$.)

 Let $G$ be a graph of girth $\ge5$, order $n=47$,
and size $6n-126 = 156$. By lemma~\ref{L:delbort}
(for $l=0$), $\mval G \ge 6$; whence $\mxval G \le
7$ by lemma~\ref{L:vv}, whence $\ds G =
(\cdots0,17,30,0\cdots)$. Now, if there were no
edges in $G_6$, then we would have $\ca{E(V_6,V_7)}
= 102$ and $e(G_7)=54$, and lemma~\ref{L:p-estim3}
(for $r=6$ and $s=3$) would yield
 $$255 = 17\cdot15 \le p(V_6,V_7) \le 435 - 54 -
12\cdot3 - 18\cdot6 = 237,$$
 a contradiction. Thus, instead, we may choose a $v
\in V_6$, with a \nb\ $w_1,$ which also has degree
6. By lemma~\ref{L:delbort} (this time with $l=1$),
thus, indeed, $G-\{v,w_1\}$ is a subgraph of order
45 and size 145, and lemma~\ref{L:sup-HS} applies.

 Now, if instead $48 \le n=n(G) \le 49$, but still
$e(G) = 6n-126$, then $G$ contains a vertex $v$ of
degree 6, and $G-\{v\}$ is an extremal graph of
order $n-1$, which `recursively' is a supergraph of
an extremal order 45 graph. Thus, and by inspecting
the few induced subgraphs of $HS$ of orders $\ge47$,
theorem~\ref{T:HSsub} is proven for all $n \ne 40$.

 \section{Bounds for \T\cle4,n , for $33 \le n \le
39$.}
 \label{S:bounds}

 \begin{table}
 \caption{\it The lower known $\T\cle4,n $ (OEIS
{\bf A006856} in October, 2015).}
 \label{tab:low}
 $$\begin{array}{|c|llllllllllllllll|}
 \hline
 n&1&2&3&4&5&6&7&8&9&10&11&12&13&14&15&16\\
 \hline
 \T\cle4,n
&0&1&2&3&5&6&8&10&12&15&16&18&21&23&26&28\\
 \hline \noalign{\vskip10pt} \hline
 n&17&18&19&20&21&22&23&24&25&26&27&28&29&30&31&32\\
 \hline
 \T\cle4,n
&31&34&38&41&44&47&50&54&57&61&65&68&72&76&80&85\\
 \hline
 \end{array}$$
 \end{table}

 \begin{table}
 \caption{\it The higher known $\T\cle4,n $.}
 \label{tab:high}
 $$\begin{array}{|c|lllllllllll|}
 \hline
 n&40&41&42&43&44&45&46&47&48&49&50\\
 \hline
 \T\cle4,n
&120&124&129&134&139&145&150&156&162&168&175\\
 \hline
 \end{array}$$
 \end{table}

 \begin{table}
 \caption{\it To my knowledge the narrowest known
intermediate $\T\cle4,n $ potential value ranges.}
 \label{tab:im}
 $$\begin{array}{|c|ccccccc|}
 \hline
 n&33&34&35&36&37&38&39\\
 \hline
 \T\cle4,n & 87-89 & 90-93 & 95-98 & 99-103 &
104-107 & 109-112 & 114-116\\
 \hline
 \end{array}$$
 \end{table}

 \bigskip

 For proving lemma~\ref{L:m2}, we shall make
recursive use of tables~\ref{tab:low}
and~\ref{tab:im}. The former is well-known; the
latter proven in this section.

 Actually, also six of the seven lower bounds are
known, since they were tabled in~\cite{GN}; the
exception is that only $\T\cle4,35 \ge 94$ was
established there. On the other hand, for all
$n\ge34$, the lower bounds given in
table~\ref{tab:im} are realised by $HS$ subgraphs;
but for $n=33$ I only found such subgraphs of
sizes~$\le86$.

 In the rest of this section, we consider the upper
bounds in table~\ref{tab:im}, and we let $G$ be a
girth $\ge5$ graph with $n(G)=n$, $e(G)=e$, $\mval G
=\mv$, and $\mxval G =\mx$. The upper bounds are
determined by increasing $n$; whence we may employ
the bounds on the \T\cle4,m \ for $m<n$ when
treating $n$. Some bounds follow by just the
standard argument that a higher $e$ or a lower $\mv$
would violate the inequalities
 $${2e\overwithdelims\lfloor\rfloor n} \ge \mv \ge
e-\T\cle4,n-1 .$$

 If $(n,e) = (33,90)$, then $5 \ge \mv \ge
e-\T\cle4,32 = 5 \implies \mv=5$. Thus, any $v \in
V_{\ge7}$ would have $\vv2v \ge 7\mv = 35 > 32$
contradicting lemma~\ref{L:vv}, whence instead
$\mx=6$. Hence, repressing all trailing zeros in
degree sequences, $\ds G = (18,15)$; and all \nb s
$w$ of a $v \in V_5$ must be of degree 6, since $\v
v+\v w \ge e(G)-\T\cle4,31 + 1 = 11$ by
lemma~\ref{L:delbort}. Thus, $G$ were bipartite and
(5,6) biregular, with parts $V_5$ and $V_6$. For a
$v \in V_5$, this would force $p_3(\{v\},V_6,V_5) =
25 > 17 = \ca{V_5\setminus\{v\}}$, a contradiction.

 \medskip
 Eliminating $(n,e) = (34,94)$ takes more effort.
Step by step, we provide more and more structural
conditions on the induced subgraphs
 $$G_i = G_i(G) := G[V_i(G)] :=G[\{v \in V(G) : \v
v=i\}],$$
 on {\sl their\/} induced subgraphs $G_j(G_i)$, on
unions of corresponding vertex sets, {\sl et
cetera}, at last deducing a contradiction. 

 Indeed, arguing as in the $n=33$ case, $\mv = 5$
and $\mx=6$, and by lemma~\ref{L:delbort} also
$\mxval G_5 \le 1$. Thus, with $e_5 := e(G_5)$ and
$z := \ca {E(V_5,V_6)}$, we have
 $$\ds G = (16,18) \land (z ,e(G_6)) =
(80-2e_5,14+e_5) \land 0 \le e_5 \le 8\,.$$
  Thus, and by lemma~\ref{L:p-estim3} applied for
$(V',V'') = (V_5,V_6)$, and since, for $e_5=8$, the
first and last expressions in
lemma~\ref{L:p-estim3}\ differ just~1, forcing an as
even as possible induced degree sequence in $G_6$,
 \begin{equation}
 \label{E:34-1}
 e_5=8 \land \ds{G_6} = (10,8) \land G_5 = G_1(G_5)
= 8\p2\,,
 \end{equation}
 the disjoint union of eight order two paths, which
we may name $\p2'$, $\p2''$\kdots $\p2^{(8)}$.

 Now, $(v \in V_3(G_6) \implies \vv2v=33=n-1
\implies V_5 \subset V = \B2v)$ by lemma~\ref{L:vv};
which, for such a $v$, since $p_3(\{v\},V_5,V_5) =
\ca{V_5\cap S(v)} = 3$, forces $p_3(\{v\},V_6,V_5) =
16-3-3 = 10$, and thus the induced second degree
$\vvind{G_6}2v = 6\ca{V_6\cap S(v)}-10 = 8$; i.~e.,
that
 $$G_3(G_6) = G_2(G_3(G_6)) = \c8$$
 (the 8-cycle graph, with $V(\c8) = \{\row c8\}$,
say).
 Similarly, for any $v \in V_2(G_6)$ there is a
single $z \in V$ at distance 3 from $v$, of degree
either 5 or~6, and with $v$ having $\vvind{G_6}2v =
5$ or~4, respectively;
whence
 $$V_2(G_6) = V_1(G_2(G_6)) \cup V_2(G_2(G_6)) \land
\ca{V_1(G_2(G_6))} = 8 \land \ca{V_2(G_2(G_6))} =
2,$$
 whence in particular
 \begin{equation}
 \label{E:34-2}
 G_2(G_6) \hbox{ has three or two $\p2$
components;}
 \end{equation}
 the amount depending on whether or not
$G_2(G_2(G_6))$ is connected.
 
 The girth property and~(\ref{E:34-1}) yield that
 \begin{equation}
 \label{E:34-3}
 \hbox{the $V_5$ \nb s of any edge in } E(V_6)\
\hbox{are in different $G_5$ components.}
 \end{equation}

 Next, we explicitly match the eight $c_i \in V(\c8)
= V_3(G_6)$ and the eight components $\p2^{(i)}$ of
$G_5\,$ (in both cases counting indices modulo 8).
Let $\{d_i\} := V_2(G_6) \cap S(c_i)$. $c_i$ is
adjacent to three $G_5$ components. Since each
degree 5 vertex is the end vertex of exactly one \p3
or \p2 from $c_i$, the $V_5$ \nb s of $d_i$,
$c_{i-1}$, and $c_{i+1}$ together form the remaining
five $G_5$ components. More precisely, $d_i$ is
adjacent to four of these, while the fifth
component, and only that component, is adjacent to
both $c_{i-1}$ and $c_{i+1}$.

 Thus, any three successive elements in \c8 must be
adjacent to $3+3+(3-1) = 8$ of the $G_5$ components,
i.~e., to all of them. Hence, if we consider four
successuve elements $c_{i-1}$, $c_i$, $c_{i+1}$, and
$c_{i+2}$, then both $c_{i-1}$ and $c_{i+2}$ must be
adjacent to both the $G_5$ components adjacent to
neither $c_i$ nor $c_{i+1}$; whence $c_{i-1}$ and
$c_{c+2}$ must have exactly these two components as
common component \nb s; the third \p2 adjacent to
$c_{i+2}$ also is adjacent to $c_i$, and thus
by~(\ref{E:34-3}) not to $c_{i-1}\,$.

 Thus, $c_i$ shares two of its three adjacent $G_5$
components with $c_{i+3}$, and two with $c_{i-3}$,
whence, by the principle of inclusion-exclusion, all
three of $c_i$ and the $c_{i\pm3}$ share one $G_5$
component \nb, say $\p2^{(\nu(i))}$. On the other
hand, by~(\ref{E:34-3}), $\p2^{(\nu(i))}$ has no
other \nb\ in \c8. Thus, $\nu$ is injective, and
hence bijective, whence without loss of generality\
we may let it be the identity. Summing up, and also
employing that $S(c_{i-3}) \cap S(c_{i+3})$ consists
of {\sl only\/} $c_{i+4}\,$, we find that
 \begin{equation}
 \label{E:34-4}
 \hbox{the } \p2^{(i)} \hbox{ \nb s in } V_3(G_6)
\hbox{ are } \{c_{i-3}, c_i, c_{i+3}\}, \hbox{ and }
\p2^{(i)} \subset S(c_{i-3}) \cup S(c_{i+3}),
 \end{equation}
 for all $i$.

 However, by~(\ref{E:34-2}), $G_{6,2}$ contains a
\p2 component, which without loss of generality\ be
$\{d_1,d_i\}$, where $i \in \{3\kdots7\}$ by the
girth condition. Each $G_5$ component is adjacent to
exactly one \p2 vertex; and by~(\ref{E:34-3})
applied for \2c_1,d_1. and for \2c_i,d_i.,
respectively, $\p2'$, $\p2^{(4)}$, and $\p2^{(6)}$
are adjacent to $d_i$, but $\p2^{(i)}$ and the
$\p2^{(i\pm3)}$ to $d_1$. Thus $\{1,4,6\} \cap
\{i-3,i,i+3\} = \0$, forcing $i=5$. Thus and without
loss of generality, $p_1^{(5)} \in \p2^{(5)} \cap
S(d_1)$; but by~(\ref{E:34-4}) $p_1^{(5)}$ also is
adjacent to a $c_j \in \{c_2, c_8\}$, causing the
existence of a 4-cycle $\{c_1, d_1, p_1^{(5)},
c_j\}$, and thus the sought contradiction.

 \medskip
 Thus and by the standard argument, $\T\cle4,35 \le
98$ and $\T\cle4,36 \le 103$.

 \smallskip
 For the next few bounds, we use the fact (from
\cite{OW}) that (6,5) cages have order 40, whence
$G$ cannot be 6-regular for $n<40$.

 If $(n,e)=(37,108)$ and $\mx=6$ (whence $\ds
G=(6,31)$), then $(v,v' \in V_5 \implies \dist(v,v')
\ge 3)$ by lemma~\ref{L:stigbort}, whence we could
make a 6-regular girth $\ge5$ realiser $G'$ out of
$G$, by adding one vertex and edges from that to all
$V_5$ vertices; yielding the contradiction $n(G') =
38 < 40$.

 Thus, instead, for such graphs we should have
$\mx=7$, but $\v v=7 \implies \ca{V_5\cap S(v)} \ge
5$; contradicting $\T\cle4,31 < 81$, by
lemma~\ref{L:delbort} with $l=5$.

 Thus, instead, $\T\cle4,37 \le 107$, whence (and by
the standard argument) $\T\cle4,38 \le 112$.

 \smallskip
 For a while and for a contradiction, assume that
$(n,e) = (39,117)$. As before, $G$ is not 6-regular,
whence, instead and by the standard arguments,
 $$\mv = 5 \hbox{ and } \mx = 7\,.$$

 Next, consider a $v \in V_7$, with \nb s \row w7,
say, where we may assume $\v{w_1} \le \ldots \le
\v{w_7}$. Since then, on the one hand, $38 \ge
\vv2v$ by lemma~\ref{L:vv}, while on the other hand
 $$\sum\3i15 \v{w_i} \ge 26$$
 by lemma~\ref{L:delbort}, then $\v{w_1} = \v{w_2} =
\v{w_3} = \v{w_4}=5$, but $\v{w_5} = \v{w_6} =
\v{w_7} = 6$.
 In particular, each $v \in V_7$ has four \nb s in
$V_5$, and all its \nb s in $V_{\le6}\,$. Hence,
$n_7 > 0 \implies n_7 \ge n_5 \ge 4$.

 Thus, in applying lemma~\ref{L:p-estim3} for
$(V',V'') := (V_7,V_5)$,
 $${n_5\choose2} \ge {4\choose2}n_7 \ge 6n_5
\implies 13 \le n_5 \le n_7\,.$$
 On the other hand, since there are $7n_7$ edges
between $V_7$ and $V_{\le6}\,$, $n_7 \le 16$; and
for each of the potential four values of $n_7$,
lemma~\ref{L:p-estim3} applied for $(V',V'') :=
(V_{\le6},V_7)$ yields a contradiction.

 \section{Proof outline for the second main lemma.}
 \label{S:ml2}

 This proof of lemma~\ref{L:m2} is fairly eclectic.
It uses both the more `independent' methods of the
proof of lemma~\ref{L:m}, and more `conventional'
methods, including references to some upper bounds
for lower orders, and case division. Where the proof
mainly re-uses already presented ideas, it just is
outlined.

 We retain the notation from section~\ref{S:ml},
{\sl mutatis mutandis}. Thus, this time, the
globally most even degree distribution is $\es =
(\cdots,0,40,0\cdots) \in \R^{40}$, and we find that
 $$\bs = \ds G = \es \iff \mv = \mval G = 6 \iff \mx
= \mxval G = 6.$$
 For the rest of the proof (and for a
contradiction), we assume the converse, i.~e., that
$\mv\le5$ and $\mx\ge7$. Actually, since $\T\cle4,31
\le 80 < e-39$, and by lemma~\ref{L:helbort}, we
must have $\mx=7$, i.~e., that
 \begin{equation}
 a = n'' = n_7 \hbox{ and } b=0\,.
 \end{equation}
 Likewise, employing $\T\cle4,39 \le 116$,
$\T\cle4,37 \le 107$, and $\T\cle4,38 \le 112$, and
putting $\tilde n = \ca{V_{\le5}}$, we get
 \begin{equation}
 \label{E:c}
 \mv\ge4,\ c := n_4 \le 2,\ e(G_4) = 0, \hbox{ and }
a = \tilde n+c\,.
 \end{equation}
 In other words,
 \begin{equation}
 \label{E:amax}
 \bs = \ds G = (\cdots0,c,a-2c,40-2a+c,a,0\cdots),
\hbox{ whence } a \le 0.5(40+c) \le 21\,.
 \end{equation}
 Moreover, if both $v$ and three of its \nb s would
have degree 7, then the sum of the degrees of the
other four \nb s of $v$ by lemma~\ref{L:vv} would
not exceed $39-3\cdot21 = 18$, which would
contradict lemma~\ref{L:delbort}, since $\T\cle4,35
\le 98 < e-7-18+4$. Thus, instead, and by counting,
 \begin{equation}
 \label{E:z}
 \mxval G_7 \le 2\implies e'' \le a \land z =
7a-2e'' \ge 5a\,.
 \end{equation}
 We shall make no further `recursive' reference to
bounds for \T\cle4,m \ for any $m<40$.

 \medskip
 This time, we shall have no use for virtually
trunkated degree sequences.

 For $p := p(V',V'') = p(V_{\le6},V_7)$, consider
the upper estimate
 $$p_u := {a\choose2} - e'' - \max(0, 2(e''-0.5a)) =
0.5\min(a^2-a-2e'', a^2+a-6e'') \ge p,$$
 the minimum of the two upper bounds in
lemma~\ref{L:p-estim3} got by choosing $s=0$ and
$s=1$, respectively.

 As our first crude lower $p$ estimate, we similarly
may combine the lower estimates of that lemma for
$r=1,2,3$, putting
 \begin{eqnarray*}
 p_l := \max(z-n', (3n'-z)\cdot1+(z-2n')\cdot3,
(4n'-z)\cdot3+(z-3n')\cdot6 )\\
  = \max (8a-2e''-40, 17a-4e''-120, 27a-6e''-240).
 \end{eqnarray*}
 In particular, and by~(\ref{E:z}), employing the
$s=1$ and the $r=3$ clauses,
 $$0 \le p_u-p_l \le 0.5(a^2+a-6e'')-(27a-6e''-240)
= 0.5(a^2-53a+6e''+480) \le 0.5 (a-15) (a-32).$$
 Thus, and by~(\ref{E:amax}),
 \begin{equation}
 \label{E:X14}
 a \le 15;
 \end{equation}
 and for $a=15$, only the maximal $e''=a$ would be
possible. Similarly, by means of the $r=2$ and both
the $s$ clauses, we get restrictions on $e''$
for $a\ge10$:
 \begin{eqnarray*}
 a \ge 10 \implies e'' &\ge {\textstyle \max \left(
\lceil \frac16 (35a-a^2-240) \rceil , \lceil \frac12
(33a-a^2-240) \rceil \right)}\\
 &= \max \left( 2a-18,
16-{17-a\choose2} \right).
 \end{eqnarray*}

 Both (\ref{E:amax}) and (\ref{E:X14}) depend only
on the crude $p_l$ we got, assuming an even
distribution of $\deg''$ on $V_{\le6}\,$. However,
in average, $\deg''$ is higher on $V_{\le5}$ than on
$V_6$, which leads to sharper bounds, by considering
$p(V_{\le5},V_7)$ and $p(V_6,V_7)$ separately, and
noting that
 $$p = p(V_{\le5},V_7) + p(V_6,V_7)$$
 by lemma~\ref{L:p-estim2}. In fact, with
$\tilde z := \ca{E(V_{\le5},V_7)}$,
 \begin{equation}
 \label{E:tld-z}
 \tilde z \ge 3a-4c+2e''.
 \end{equation}
 This enables us to work with improved higher $p_l$,
and (by worst case analyses for the various $a$,
combined with analysis of the impact of rising $e''$
above the~(\ref{E:tld-z}) bound), we may sharpen
(\ref{E:X14}) to
 $$a \le 8.$$
 Moreover, the same analysis yields
 $$a=8 \implies c=2 \land e''=0 \land \max\limits_{v
\in V_{\le5}} \deg''(v) = 3 \land \tilde z = 16.$$
 However, this $\tilde z$ value is the minimal one
allowed in~(\ref{E:tld-z}), and a closer analysis of
that inequality reveals that the term $-4c$ therein
may be replaced by
 $$-\sum_{v\in V_4} \deg''(v);$$
 a quantity that thus on the one hand were at least
$-6$, and on the other hand were equal to $-8$; a
contradiction. Thus, instead, and summing up,
 \begin{equation}
 \label{E:amax2}
 a \le 7.
 \end{equation}

 For the lowest $a$ values, also work with $\tilde p
:= p(V_7,V_{\le5})$. In order to force a larger
$\tilde z$ than given by~(\ref{E:tld-z}), we
sometimes may employ a {\sl discharging
technique\/}: We start by giving each element in
$V_5\choose2$ (i.~e., 2-subset in $V_5$) a charge
$\frac13$, and every other element in
$V_{\le5}\choose2$ a charge~1. Next, we move the
charge of any 2-subset of $V_{\le5}$ with a common
\nb\ $v \in V_7$ to that $v$. Then, after
discharging, each vertex of degree 7 has received a
total charge $\ge1$, since $v$ either has precisely
two \nb s in $V_{\le5}\,$, of which at least one has
degree 4, or has at least three \nb s in $V_5\,$.
Moreover, the number of 3-subsets of $n_5$ belonging
to the \nb hoods of different $v_7$ elements is at
most
 $$\left\lfloor{\textstyle{n_5\over3} \cdot
\lfloor0.5(n_5-1)\rfloor} \right\rfloor - \ep,
\hbox{ where } \ep = \cases {1&if
$n_5\equiv5\pmod6$\cr 0&else\cr},$$
 since different such 3-sets share at most one
vertex, and by~\cite[theorem~1]G.

 This eliminates almost all the remaining
possibilities. The few exceptions are treated by
structure determination and case division,
eventually leading to dismissal. The most complex of
these treatments occur for the cases where $a=7$,
$c=2$, and $e''=0$; let us briefly consider them. In
these cases, we may put $V_7 = \{\row v7\}$, $V_4 =
\{u_1,u_2\}$, $V_5 = \{w_1,w_2,w_3\}$, and $Z_i =
V_{\le5} \cap S(v_i)$ for $i=1\kdots7$. Note, that
then $\tilde z \ge 14$, since each $\ca{S(v_i)} \ge
2$. On the other hand, $\tilde z \le 15$, since
$p(V_7,V_{\le5}) \le {5\choose2} < 5\cdot1+2\cdot3$.
This yields two cases: $\tilde z = 14$ and $\tilde z
= 15$. In either case,
 $$p(V_6,V_7) \ge z-\tilde z-n_6 = 21-\tilde z
\implies p(V_{\le5},V_7) \le \tilde z.$$

 However, if here $\tilde z = 14$, then each one of
the seven elements in $V_7$ must be adjacent to
$V_4$, whence and without loss of generality $S(u_1)
\subset V_7$, whence each one of the four 2-subsets
of $V_{\le5}$ containing $u$ were employed as some
$Z_i$, whence in particular some $Z_i =
\{u_1,u_2\}$, whence $S(u_2) \subset V_7$, too;
yielding
 $$p(V_{\le5},V_7) = 3\cdot1+2\cdot6 = 15 > 14 =
\tilde z,$$
 a contradiction.

 In the remaining case, $\tilde z = 15$; any $S(u_i)
\subset V_7$ again would yield a too high
$p(V_{\le5},V_7)$, and contradiction; whence instead
(and without loss of generality) $Z_7 = V_5$, and
 $$Z_i = \{u_j,w_k\}\ \hbox{ with }\ i \equiv j
\pmod2\ \hbox{ and }\ i \equiv k \pmod3$$
 for $j=1\kdots6$. By the girth condition
and~(\ref{E:c}),
 $$E(V_{\le5}) = E(V_4) = \0.$$
 Thus, we have determined the graph structure of
$G(V_{\le5} \cup V_7)$ and may continue to the
degree 6 vertices. Put
 $$A := V_6 \cap \bigcup\3i13 S(w_i),\ B = \{\row
b4\} := V_6 \cap S(v_7),\ C := V_6 \setminus(A \cup
B),$$
 and put $C_i := V_6 \cap S(b_i)$ for $i=1\kdots4$.
 By the girth conditions, $A \cap B = E(A,B) = \0$
and $\ca A = 3\cdot2$;
whence $(A,B,C)$ is a tripartition of $V_6$, with
 $$(\ca A, \ca B, \ca C) = (4,6,18).$$
 On the other hand, the $b_i$ together have 20 \nb s
apart from $v_7$, and all these \nb s are different.
Hence,
 $$\bigcup\3i14 S(b_i) = C \cup V_4 \cup \{v_7\}.$$
 Thus, $\{C_1,C_2,C_3,C_4\}$ is a partition of $C$
into four parts, of sizes $5+5+4+4$ or $5+5+5+3$.

 We now continue with the degree 6 \nb s of the
$v_i$, for $1 \le i \le 6$. For these $i$,
$\ca{V_6\cap S(v_i)} = 5$. However, $B \cap S(v_i) =
\0$ (girth reasons), and $v_i$ has at most one \nb\
in each $C_j$. Moreover, some $u_l \in S(v_i)$, and
some $b_{j'} \in S(u_k)$, and $v_i$ has no \nb\ in
$C_{j'}$ (girth reasons). Thus, $\ca {C\cap S(v_i)}
\le 3$; whence $\ca {A\cap S(i)} \ge 2$. Summing up,
and by counting,
 $$\ca {B\cap C} = 18 \land E(V_{\le5},C)
= \0 \land \ca {E(V_7,C)} \le 18 \land \ca {E(A,C)}
\le 18 \implies e(C) \ge 27.$$
 On the other hand, no $C$ vertex could have a \nb\
in its own part; and for $1 \le i < j \le 4$, $\ca
{E(C_i,C_j)} \le \min (\ca{C_i}, \ca{C_j})$, whence
 $$e(C) \le \max(1\cdot5+5\cdot4, 3\cdot5+3\cdot3) =
25 < 27,$$
 the sought contradiction.

 \medskip
 After similarly having eliminated all remaining
cases with $\mx=7$, thus, indeed, only the
possibility that $G$ is 6-regular remains.\noproof

 In particular and by \cite W we have
 \begin{corollary}
 \label{C:HSsub_40}
 The graphs treated in lemma~\ref{L:m2} are unique
(up to isomorphisms), and are subgraphs of
Hoffman-Singleton graphs.\noproof
 \end{corollary}

 In particular, this concludes the proof of
teorem~\ref{T:HSsub}.

 \section{Proof of theorem~\ref{T2}.}
 \label{S:thp2}

 From lemma~\ref{L:m2}, we indeed may deduce
theorem~\ref{T2} in a few steps, which we briefly
indicate. In each case, we assume that $G$ has girth
$\ge5$, order~$n$, and size one more that the
proposed value of \T\cle4,n . We sometimes refer to
the lemma~\ref{L:p-estim3} upper and lower bounds of
$p := p(V_{\le6},V_7)$ as $p_u$ and $p_l$,
respectively.

 Start by noting that indeed $\T\cle4,40 < 121$;
proven as for \T\cle4,45 . Thus:

 $n=41\ ($and $e(G)=124+1) \implies \mx\ge7 \implies
\mv=5 \land n_5\ge2$; but then $v \in V_5 \implies
G-\{v\}$ not 6-regular, contradicting
lemma~\ref{L:m2}. Thus:

 $n=42 \implies \mx\ge7$, contradicting $ \mv>5$.
Thus:

 $n=43 \implies \mx=7 \land G_7 \simeq
\overline{K_{12}}$ (the edge-free graph of order
12), yielding $p_u-p_l \le 66 - (9\cdot1 + 22\cdot3)
= -9 < 0$, which contradicts $p_l \le p \le p_u\,$.
Thus, finally:

 $n=44 \implies \mx=7 \land n_7=16 \land e(G_7) \le
8 \implies p_u-p_l \le (120-8) - (16\cdot3+12\cdot6)
= -8 < 0$, a contradiction.\noproof

 \bibliographystyle{plain}

 \end{document}